\documentclass[11pt, leqno,twoside]{article}
\usepackage{amssymb}
\usepackage{amsmath}
\usepackage{amsthm}
\usepackage{amsfonts}
\usepackage{color}
\usepackage{textcomp}
\usepackage{graphicx}

\usepackage[active]{srcltx}

\allowdisplaybreaks

\def\XXint#1#2#3{{\setbox0=\hbox{$#1{#2#3}{\int}$ }
\vcenter{\hbox{$#2#3$ }}\kern-.6\wd0}}

\begin{document}

\title{Uniqueness in Calder\'{o}n Problem with Partial Data for Less Smooth Conductivities
\thanks
{E-mail addresses: \   guozhang@jyu.fi(G. Zhang)}
}
\author{{ $\text{ Guo Zhang}$} \\
{\small Department of Mathematics and Statistics, P.O. Box 35 (MaD), FI-40014}\\
{\small University of Jyv\"{a}skyl\"{a}, Finland }}

\date{}
\maketitle

\begin{abstract}
In this paper we study the inverse conductivity problem with partial data. Moreover,  we show that, in dimension $n\geq 3$  the uniqueness of the  Calder\'{o}n
problem holds  for the $C^{1}\bigcap H^{\frac{3}{2},2}$ conductivities.

 \end{abstract}
{\bf MR Subject Classification:}\ \ 35R30.

\bigbreak
\begin{center}
\textbf{\Large {\S 1.\ Introduction}}
\end{center}

In 1980, Calder\'{o}n \cite{Cal} considered whether one can determine the electrical conductivity of a medium by making voltage and current measurements at the boundary 
of the medium. This inverse method is known as Electrical impedance tomography (EIT). EIT also arises in medical imaging given that human organs and tissues have quite
 different conductivities \cite{Jos}. One exciting potential application is the early diagnosis of breast cancer \cite{ZG}. The conductivity of a malignant tumor is typically
$0.2$ mho which is significantly higher than normal tissue which has been typically measured at $0.03$ mho. Another application is to monitor pulmonary function \cite{INGC}.
See the book \cite{Hol} and the issue of Physiological measurement \cite{IMS} for other medical imaging application of EIT.

 We now describe more precisely the mathematical problem.

 Let $\Omega\subset\mathbb{R}^n, n\geq3$ be an open, bounded domain with $C^2$ boundary $\partial\Omega$, and let $\gamma$ be a strictly positive real valued function
 defined on $\Omega$ which gives the conductivity at a given point. Given a Voltage potential $f$ on the boundary, the equation for the potential in the interior, under
the assumption of no sinks or sources of current in $\Omega$, is
\[
\renewcommand{\arraystretch}{1.25}
\begin{array}{lll}
\text{div}(\gamma\nabla u)=0,\ \ \ \ \text{in}\ \Omega,\ \ \ \ u$\textbar$_{\partial\Omega}=f.
\end{array}
\eqno(1.1)
\]
The Dirichlet-to-Neumann map is defined in this case as follows:
\[
\renewcommand{\arraystretch}{1.25}
\begin{array}{lll}
\Lambda_{\gamma}(f)=\gamma\frac{\partial u}{\partial \nu}$\textbar$_{\partial\Omega},
\end{array}
\eqno(1.2)
\]
where $\frac{\partial}{\partial \nu}$ is the outward normal derivatives at the boundary. For $\gamma\in \text{Lip}(\bar\Omega)$, then
 $\Lambda_{\gamma}$ is a well defined map
from $H^{\frac{1}{2}}(\partial\Omega)$ to $H^{-\frac{1}{2}}(\partial\Omega)$.

The Calder\'{o}n problem concerns the inversion of the map $\gamma\rightarrow\Lambda_{\gamma}$, i.e., whether $\Lambda_{\gamma}$ determines
 $\gamma$ uniquely and in that
case how to reconstruct $\gamma$ from $\Lambda_{\gamma}$.

For the uniqueness, it was first proved for smooth conductivities by Sylvester and Uhlmann in their fundamental paper \cite{SU1}, which opened the door
of studying the Calder\'{o}n problem.

Now we briefly recall the basic idea in \cite{SU1}. For the $C^2$ conductivity $\gamma$, letting $u=\gamma^{\frac{1}{2}}v$ be  the solution of $(1.1)$, we can
 deduce $v$ satisfies the Schr\"{o}dinger equation
\[
\renewcommand{\arraystretch}{1.25}
\begin{array}{lll}
(-\triangle+q)v=0,
\end{array}
\eqno(1.3)
\]
where $q$ is defined by $q=\gamma^{-\frac{1}{2}}\triangle\gamma^{\frac{1}{2}}$. The corresponding Dirichlet-to-Neumann map is defined by
 $\Lambda_{q}(f)=\frac{\partial v}{\partial \nu}$ with the boundary data $f$. To be more important, for the $C^2$ conductivities
 $\gamma_1$, $\gamma_2$, $\Lambda_{\gamma_1}=\Lambda_{\gamma_2}$ implies $\Lambda_{q_1}=\Lambda_{q_2}$ for
 $q_i=\gamma_{i}^{-\frac{1}{2}}\triangle\gamma_{i}^{\frac{1}{2}}, i=1,2$ (see \cite{KV1} \cite{A1}).  Therefore,  we convert the inverse problem for the conductivity
 equation $(1.1)$ to an equivalent inverse problem for the Schr\"{o}dinger equation $(1.3)$.

 If $\Lambda_{q_1}=\Lambda_{q_2}$ and $(-\triangle+q_i)v_i=0$, then a simple calculation shows that
 \[
\renewcommand{\arraystretch}{1.25}
\begin{array}{lll}
\displaystyle\int_{\Omega}(q_1-q_2)v_1v_2dx=0.
\end{array}
\eqno(1.4)
\]
From this discussion,  to get $q_1=q_2$, we just need construct many enough solutions to the corresponding  Schr\"{o}dinger equations $(1.3)$ such that
their products are dense in some sense.

In \cite{SU1}, they construct this type of complex geometrical optics solutions $v_i=e^{x\cdot\zeta_i}(1+w_i)$, where the $\zeta_i \in C^n$ are
 chosen so that $\zeta_i\cdot\zeta_i=0$, which implies $e^{x\cdot\zeta_i}$ is harmonic and $e^{x\cdot\zeta_1}e^{x\cdot\zeta_2}=e^{i x\cdot k}$
 for some fixed frequency $k\in \mathbb{R}^n$. For $w_i$, these are the solutions of the following equations
 \[
\renewcommand{\arraystretch}{1.25}
\begin{array}{lll}
\triangle_{\zeta_i}w_i:=\triangle w_i+2\zeta_i\cdot\nabla w_i=q_i(1+w_i).
\end{array}
\eqno(1.5)
\]

In three or more dimensions, we can find infinite high frequency solutions $v_i=e^{x\cdot\zeta_i}(1+w_i)$ satisfying the remainders $w_i$ decay to
zero in some sense as $|\zeta_i|\rightarrow \infty$, so that the product $v_1v_2$ converges to $e^{ix \cdot k}$. Uniqueness then follows from Fourier  inversion.

For the less smooth conductivities $\gamma_i$, following this idea, one may find the solutions $w_i$ of $(1.5)$ in some suitable  Sobolev or
 Besov space by contract mapping theorem. Furthermore,  the chosen CGO solution $v_i$ such that $(1.4)$ makes sense. In view of the above analysis, one may show the
 uniqueness of the Calder\'{o}n problem  which holds under different types of conductivities. For example, Brown \cite{B1} obtain uniqueness under the assumption
of $\frac{3}{2}+\epsilon$
derivatives. Later,  uniqueness for exactly $\frac{3}{2}$ bounded derivatives was shown in \cite{PPU1} and for $\frac{3}{2}$ derivatives being in $L^p, p>2n$
was shown in \cite{BT1}.

In another direction, for some special conductivity, Greenleaf,  Lassas, and Uhlmann \cite{GLU1} obtained global uniqueness for certain
conductivities in $C^{1+\epsilon}$. Later,   Kim \cite{K2} established global uniqueness for Lipschitz conductivities that are piecewise smooth across
polyhedral boundaries.

Following the above idea, the sharpest uniqueness results so far seem to require essentially $\frac{3}{2}$ derivatives of the conductivity.
 Recently, Haberman and Tataru \cite{HT1}  use a totally new idea to show uniqueness for almost Lipschitz conductivities. The main idea is as follows.
Since the symbol of the operator $\triangle_\zeta$ is $-|\zeta|^2+2i\zeta\cdot\xi$,  Haberman and Tataru introduce Bourgain's  spaces 
$\dot{X}^{s,\frac{1}{2}}$ \cite{Bou}, which are defined by the norm $\|u\|_{\dot{X}^{s,\frac{1}{2}}_{\zeta}}=\||P_{\zeta}(\xi)|^{\frac{1}{2}}\hat u(\xi)\|_{L^2}$,
where $P_{\zeta}(\xi)=-|\zeta|^2+2i\zeta\cdot\xi, |\zeta|=s$. Furthermore, by contract mapping theorem, one may find a family of high frequent CGO solutions
$v_i=e^{x\cdot\zeta_i}(1+w_i)$ satisfying $\|w_i\|_{\dot{X}^{s,\frac{1}{2}}_{\zeta_i}}\rightarrow 0$, as $|\zeta_i|\rightarrow\infty$. The most important thing
is that $(1.4)$ makes sense. Uniqueness then follows from Fourier  inversion.

We will use this idea to show that the uniqueness of the Calder\'{o}n
problem  with partial data for the $C^{1}(\bar{\Omega}) \bigcap H^{\frac{3}{2},2}(\Omega)$ conductivities in this paper.

Before stating the theorem, let us recall what is known about the  uniqueness of the Calder\'{o}n
problem  with partial data.

First we introduce some notations. Let $\eta\in S^{n-1}$ and define the subsets
 \[
\renewcommand{\arraystretch}{1.25}
\begin{array}{lll}
\partial\Omega_{+}=\{x\in\partial\Omega,\ \nu(x)\cdot \eta \geq 0 \},\ \ \ \ \ \
\partial\Omega_{-}=\{x\in\partial\Omega,\ \nu(x)\cdot \eta < 0 \},
\end{array}
\]
where $\nu(x)$ is the outward normal direction at the boundary point $x$. For $\epsilon>0$, define further subsets
 \[
\renewcommand{\arraystretch}{1.25}
\begin{array}{lll}
\partial\Omega_{+,\epsilon}=\{x\in\partial\Omega,\ \nu(x)\cdot \eta \geq \epsilon \},\ \ \ \ \ \
\partial\Omega_{-,\epsilon}=\{x\in\partial\Omega,\ \nu(x)\cdot \eta < \epsilon \}.
\end{array}
\]

Bukhgeim and Uhlmann \cite{BU1} first established the uniqueness of the Calder\'{o}n
problem under assumptions $\gamma_1, \gamma_2\in C^2(\bar\Omega)$, $\gamma_1
$\textbar$_{\partial\Omega}= \gamma_2
$\textbar$_{\partial\Omega}$, and $\Lambda_{\gamma_1}$\textbar$_{\partial\Omega_{-,\epsilon}}=\Lambda_{\gamma_2}$\textbar$_{\partial\Omega_{-,\epsilon}}$.
Later Knudsen \cite{K1} generalized their result and established the uniqueness of the Calder\'{o}n problem under the
assumptions $\gamma_1, \gamma_2\in W^{\frac{3}{2}+r,2n}(\Omega)$, $r>0$, $\gamma_1
$\textbar$_{\partial\Omega_{+}}= \gamma_2
$\textbar$_{\partial\Omega_{+}}$, $\partial_{\nu}\gamma_1
$\textbar$_{\partial\Omega_{+}}= \partial_{\nu}\gamma_2
$\textbar$_{\partial\Omega_{+}}$ and  $\Lambda_{\gamma_1}$\textbar$_{\partial\Omega_{-,\epsilon}}=\Lambda_{\gamma_2}$\textbar$_{\partial\Omega_{-,\epsilon}}$.

Now we state our result as follow.

\textbf{Theorem $1.1$}\ \ Let $\Omega\subset\mathbb{R}^n, n\geq3$, be an open, bounded domain with $C^2$ boundary. For $i
=1, 2$, let $\gamma_i \in C^{1}(\bar{\Omega}) \bigcap H^{\frac{3}{2},2}(\Omega)$ be real valued function and $\gamma_i>c>0$. Fix $\eta\in S^{n-1}$ and
suppose further that $\gamma_1
$\textbar$_{\partial\Omega_{+}}= \gamma_2
$\textbar$_{\partial\Omega_{+}}$, $\partial_{\nu}\gamma_1
$\textbar$_{\partial\Omega_{+}}= \partial_{\nu}\gamma_2
$\textbar$_{\partial\Omega_{+}}$, and for some
 $\epsilon>0$,
 \[
\renewcommand{\arraystretch}{1.25}
\begin{array}{lll}
\Lambda_{\gamma_1}$\textbar$_{\partial\Omega_{-,\epsilon}}=\Lambda_{\gamma_2}$\textbar$_{\partial\Omega_{-,\epsilon}},\ \ \ \text{ for any}
 f \in H^{\frac{1}{2}}(\partial\Omega).
\end{array}
\]
Then we have $\gamma_1=\gamma_2$.

For the above known results about the Calder\'{o}n problem with partial data, the main idea is that one can use the linear  limiting   Carleman
 weight $x\cdot \eta$ in the proof.
Indeed, Kenig, Sj{\"o}strand and Uhlmann use a nonlinear limiting  Carleman  weight to obtain a new type of result in  \cite{KSU1}. By combing our
approach and  Kenig, Sj{\"o}strand and
 Uhlmann \cite{KSU1} 's idea, it is expected a generalization
of the Theorem $1.1$ can be found.

Our paper is organized as follows: In section $2$ we will state the results of approximation of conductivities. In section $3$  we will  construct the CGO solutions.
In section $4$ we will give the Carleman estimate for the CGO solutions. In section $5$ We will present the proof of Theorem $1.1$.

\bigbreak
\begin{center}
\textbf{\Large {\S 2.\   Preliminary Results }}
\end{center}

From the assumption of the conductivities $\gamma_i, i=1,2.$,  in the Theorem $1.1$, we can extend $\gamma_i$ to be the functions in the
whole space $\mathbb{R}^n$ such that
$(i)\  \gamma_i(x)\geq C$, $\ (ii)\  \gamma_i-1\in H^{\frac{3}{2},2}(\mathbb{R}^n) \ \text{with compact support} $,
$(iii)\  \gamma_i\in C^1(\mathbb{R}^n) $,  $(iv)\ \gamma_1=\gamma_2 \ \text{outside of}\  \Omega $. For the proof of this extension, the readers are recommenced
 to read the proof of Theorem $5.7$ in \cite{Sh}. Although Our assumptions on conductivities are a little bit different,  the way of the proof still works.

Let $\Psi\in C_{0}^{\infty}(\mathbb{R}^n)$ be a nonnegative radial function with $\int_{\mathbb{R}^n }\Psi dx=1, \text{spt}\Psi \subset B(0,1)$
and $\Psi_t(x)=t^n\Psi(tx)$.
Then for $\phi=\text{log}\gamma$, $A=\nabla\text{log}\gamma$,  we define $\phi_t=\Psi_t\ast\phi$, $A_t=\Psi_t\ast A$. Clearly. $A_t=\nabla \phi_t$.

Now we state some approximation results. Some of these results are taken from \cite{S1} and \cite{K1} and some are new. For the new results, we will give
the details of the proof.

\textbf{Lemma $2.1$ }\ \ Suppose $\gamma\in  C^1(\mathbb{R}^n)$ and $\gamma-1\in H^{\frac{3}{2},2}(\mathbb{R}^n)$ with compact support. Then
as $t\rightarrow\infty$, we have
\[
\renewcommand{\arraystretch}{1.25}
\begin{array}{lll}
\|{\phi_t}\|_{L^{\infty}} \leq \|\phi\|_{L^\infty}, \ \ \ \
\|A_t\|_{L^\infty} \leq \|A\|_{L^\infty},\ \ \ \ \|D^2\phi_t\|_{L^\infty}=o(t),\\
 \|\phi_t-\phi\|_{L^\infty}=o(\frac{1}{t}),\ \ \|A_t-A\|_{L^\infty}=o(1)
\end{array}
\]
and
\[
\renewcommand{\arraystretch}{1.25}
\begin{array}{lll}
\|\phi_t-\phi\|_{L^{2}(\mathbb{R}^n)}=o(\frac{1}{t^{\frac{3}{2}}}), \ \ \ \
\|\phi_t-\phi\|_{H^{1,2}(\mathbb{R}^n)}=o(\frac{1}{t^{\frac{1}{2}}}),\\
\|\phi_t-\phi\|_{H^{\frac{3}{2},2}(\mathbb{R}^n)}=o(1),\ \ \ \  \|D^{\alpha}\phi_t\|_{L^2}=o(t^{\frac{1}{2}}), |\alpha|=2.
\end{array}
\]

\textbf{Proof}\ \  Since $\gamma\in C^1(\mathbb{R}^n)$ and $\gamma-1\in H^{\frac{3}{2},2}(\mathbb{R}^n)$ with compact support,
then $\phi\in C_{c}^{1}(\mathbb{R}^n)\bigcap
H^{\frac{3}{2},2}(\mathbb{R}^n)$ and $A\in C_{c}(\mathbb{R}^n)\bigcap H^{\frac{1}{2},2}(\mathbb{R}^n)$. From Lemma $2.1$ in
chapter $2$ of  \cite{S1}, we can get
$\|{\phi_t}\|_{L^{\infty}} \leq \|\phi\|_{L^\infty}$, $\|A_t\|_{L^\infty} \leq \|A\|_{L^\infty}$, $\|D^2\phi_t\|_{L^\infty}=o(t)$
and  $\|A_t-A\|_{L^\infty}=o(1)$.

To obtain $\|\phi_t-\phi\|_{L^\infty}=o(\frac{1}{t})$, by fundamental theorem of calculus we have
\[
\renewcommand{\arraystretch}{1.25}
\begin{array}{lll}
t\biggl(\phi_t(x)-\phi(x)\biggr)&=t\biggl(\displaystyle\int_{\mathbb{R}^n}\Psi_t(x-y)\phi(y)dy-\phi(x)\biggr)\\
&=t\biggl(\displaystyle\int_{\mathbb{R}^n}\Psi(z)(\phi(x-\frac{z}{t})-\phi(x))dz\biggr)\\
&=\displaystyle\int_{\mathbb{R}^n}\Psi(z)\displaystyle\int_{0}^{1}\nabla\phi(x-s\frac{z}{t})\cdot-zdsdz.
\end{array}
\]

On the other hand, since $\Psi(z)$ is a radial function, we know
\[
\renewcommand{\arraystretch}{1.25}
\begin{array}{lll}
\displaystyle\int_{\mathbb{R}^n}\Psi(z)\nabla\phi(x)\cdot zdz=0.
\end{array}
\]

Hence,

\[
\renewcommand{\arraystretch}{1.25}
\begin{array}{lll}
t\biggl(\phi_t(x)-\phi(x)\biggr)=\displaystyle\int_{\mathbb{R}^n}\Psi(z)\displaystyle\int_{0}^{1}(\nabla\phi(x-s\frac{z}{t})-\nabla\phi(x))\cdot-zdsdz.
\end{array}
\]

Since $\text{spt}\Psi(z)\subset B(0,1)$ and $\nabla\phi\in C_{c}(\mathbb{R}^n)$, then uniform continuity  gives
\[
\renewcommand{\arraystretch}{1.25}
\begin{array}{lll}
\|\phi_t-\phi\|_{L^\infty}=o(\frac{1}{t}).
\end{array}
\]

For any $0\leq s\leq \frac{3}{2}$, by Fourier transform formula we have
\[
\renewcommand{\arraystretch}{1.25}
\begin{array}{lll}
(\phi_t-\phi)^{\hat{}}(\xi)=(\hat{\Psi}(\frac{\xi}{t})-1)\hat{\phi}(\xi).
\end{array}
\]

Hence,
\[
\renewcommand{\arraystretch}{1.25}
\begin{array}{lll}
t^{\frac{3}{2}-s}\|(\hat{\Psi}(\frac{\xi}{t})-1)|\xi|^s\hat{\phi}(\xi)\|_{L^2}=\|g(\frac{\xi}{t})|\xi|^{\frac{3}{2}}\hat{\phi}(\xi)\|_{L^2},
\end{array}
\]
where $g(z)=\frac{1}{|z|^{\frac{3}{2}-s}}(\hat{\Psi}(z)-1)$.

Since $\int_{\mathbb{R}^n}\Psi(z)dz=1$ and $\Psi$ is radial, using Fourier transform formula we obtain
\[
\renewcommand{\arraystretch}{1.25}
\begin{array}{lll}
\hat{\Psi}(0)=1,\ \ \ \ \nabla \hat{\Psi}(0)=0,
\end{array}
\]
 which implies $g(z)$ is continuous and bounded with $g(0)=0$.

Thus, Lebesgue dominated convergence theorem gives
\[
\renewcommand{\arraystretch}{1.25}
\begin{array}{lll}
\lim\limits_{t\rightarrow \infty}\|g(\frac{\xi}{t})|\xi|^{\frac{3}{2}}\hat{\phi}(\xi)\|_{L^2}=0.
\end{array}
\]

Then we obtain
\[
\renewcommand{\arraystretch}{1.25}
\begin{array}{lll}
\|\phi_t-\phi\|_{L^{2}(\mathbb{R}^n)}=o(\frac{1}{t^{\frac{3}{2}}}), \ \
\|\phi_t-\phi\|_{H^{1,2}(\mathbb{R}^n)}=o(\frac{1}{t^{\frac{1}{2}}}),\ \
\|\phi_t-\phi\|_{H^{\frac{3}{2},2}(\mathbb{R}^n)}=o(1).
\end{array}
\]

To prove $\|D^{\alpha}\phi_t\|_{L^2}=o(t^{\frac{1}{2}}), |\alpha|=2$, observe
\[
\renewcommand{\arraystretch}{1.25}
\begin{array}{lll}
t^{-\frac{1}{2}}\|D^{\alpha}\phi_t\|_{L^2}&=t^{-\frac{1}{2}}\|\xi^\alpha\hat{\Psi}(\frac{\xi}{t})\hat{\phi}(\xi)\|_{L^2}\\
&\leq \|\tilde{g}(\frac{\xi}{t})|\xi|^{\frac{3}{2}}\hat{\phi}(\xi)\|_{L^2},
\end{array}
\]
where $\tilde{g}(z)=|z|^{\frac{1}{2}}\hat{\Psi}(z)$ is continuous and bounded with $\tilde{g}(0)=0$.

Then, Lebesgue dominated convergence theorem gives that $\|D^{\alpha}\phi_t\|_{L^2}=o(t^{\frac{1}{2}}), |\alpha|=2$.
Then Lemma $2.1$ follows.

Finally, we need the following type of trace inequality for $W^{1,2}(\Omega)$. Since we can't easily find this result in the literature, 
for completeness, we will give the proof of the following Lemma.

\textbf{Lemma $2.2$}\ \ Let $\Omega\subset\mathbb{R}^n, n \geq 2$ be a $C^2$ smooth bounded domain and $u\in W^{1,2}(\Omega)$, Then there exists a constant $C$ such that 
\[
\renewcommand{\arraystretch}{1.25}
\begin{array}{lll}
\displaystyle\int_{\partial \Omega}u^2ds \leq C\biggl \{ \biggl(\displaystyle\int_{\Omega}u^2dx\biggr)^{\frac{1}{2}}\biggl(\displaystyle\int_{\Omega}|\nabla u|^2dx\biggr)^{\frac{1}{2}}
+\displaystyle\int_{\Omega} u^2dx \biggr\},
\end{array}
\]
where $C$ depend  only on $\Omega$, $n$.

\textbf{Proof}\ \ Since $\partial \Omega$ is a $C^2$ smooth boundary, we know the unit outward normal direction $\nu(x)$ is a $C^1$ vector on the boundary. Then, we
 can find a $C^1$ vector $\rho$ in $\Omega$ which is an extension of $\nu(x)$.

Thus, by divergence theorem we have 
\[
\renewcommand{\arraystretch}{1.25}
\begin{array}{lll}
\displaystyle\int_{\partial \Omega}u^2ds=\displaystyle\int_{ \Omega}\text{div}(\rho u^2) dx 
=\displaystyle\int_{ \Omega}\text{div}(\rho) u^2 dx+\displaystyle\int_{ \Omega}\rho\cdot \nabla(u^2) dx,
\end{array}
\]
in view of H\"{o}lder inequality, which implies 
\[
\renewcommand{\arraystretch}{1.25}
\begin{array}{lll}
\displaystyle\int_{\partial \Omega}u^2ds \leq C\biggl \{ \biggl(\displaystyle\int_{\Omega}u^2dx\biggr)^{\frac{1}{2}}\biggl(\displaystyle\int_{\Omega}|\nabla u|^2dx\biggr)^{\frac{1}{2}}
+\displaystyle\int_{\Omega} u^2dx \biggr\}.
\end{array}
\]

\bigbreak
\begin{center}
\textbf{\Large {\S 3. CGO Solutions}}
\end{center}

Before stating the main result, let us introduce some preliminarily results.

\textbf{Proposition $3.1$}\ \ \cite{HT1} Let $v$ and $w$ be nonnegative weights defined on $\mathbb{R}^n$. If $\Phi$ is a fixed rapidly decreasing function, then
\[
\renewcommand{\arraystretch}{1.25}
\begin{array}{lll}
\|\Phi \ast f\|_{L^2_{w}}\leq \text{min}\{\sup\limits_{\xi}\sqrt{\int J(\xi,\eta)d\eta},\  \sup\limits_{\eta}\sqrt{\int J(\xi,\eta)d\xi}\}\|f\|_{L^2_{v}},
\end{array}
\]
where
\[
\renewcommand{\arraystretch}{1.25}
\begin{array}{lll}
J(\xi,\eta)=|\Phi(\xi-\eta)|\frac{w(\xi)}{v(\eta)}.
\end{array}
\]

Now we introduce Bourgain's spaces \cite{Bou}. Using Haberman and Tataru's idea in \cite{HT1}, we define the spaces $\dot{X}^{b}_{\zeta}$ by the norm
\[
\renewcommand{\arraystretch}{1.25}
\begin{array}{lll}
\|u\|_{\dot{X}^{b}_{\zeta}}=\||P_{\zeta}(\xi)|^{b}\hat{u}(\xi)\|_{L^2},
\end{array}
\]
where $P_{\zeta}(\xi)=-|\zeta|^2+2i\zeta\cdot\xi$ is the symbol of $\triangle_{\zeta}$. In our paper, we will need the spaces $\dot{X}^{\frac{1}{2}}_{\zeta}$ and
$\dot{X}^{-\frac{1}{2}}_{\zeta}$ and also make use of inhomogeneous spaces ${X}^{b}_{\zeta}$ with norm
\[
\renewcommand{\arraystretch}{1.25}
\begin{array}{lll}
\|u\|_{{X}^{b}_{\zeta}}=\|(|\zeta|+|P_{\zeta}(\xi)|)^{b}\hat{u}(\xi)\|_{L^2}.
\end{array}
\]

Let $\zeta \in C^{n}$ be such that $\zeta\cdot\zeta=0$ and write $\zeta=s(e_1-ie_2)$, with $e_1, e_2 \in \mathbb{R}^n$ satisfying $e_1 \cdot e_2=0$.
Taking the open balls $B_0$, $B$ s.t. $\Omega\Subset B_0 \Subset B$ and choosing a schwartz cutoff function $\Phi_{B}$ which is equal to one on an open ball $B$,
then we have the following known results which are taken from \cite{HT1}.

\textbf{Proposition $3.2$}\ \ Let $\Phi_{B}$ be a fixed Schwartz function defined as  above, and write $u_{B}= \Phi_{B}u$. Then the following estimates hold with constants depending
 on $\Phi_{B}$
\[
\renewcommand{\arraystretch}{1.25}
\begin{array}{lll}
\|u_{B}\|_{\dot{X}^{-\frac{1}{2}}_{\zeta}}\lesssim  \|u\|_{{X}^{-\frac{1}{2}}_{\zeta}},
\end{array}
\eqno{(3.1)}
\]
\[
\renewcommand{\arraystretch}{1.25}
\begin{array}{lll}
\|u_{B}\|_{{X}^{\frac{1}{2}}_{\zeta}}\lesssim  \|u\|_{\dot{X}^{\frac{1}{2}}_{\zeta}},
\end{array}
\eqno{(3.2)}
\]
\[
\renewcommand{\arraystretch}{1.25}
\begin{array}{lll}
\|u_{B}\|_{L^2}\lesssim s^{-\frac{1}{2}} \|u\|_{\dot{X}^{\frac{1}{2}}_{\zeta}}.
\end{array}
\eqno{(3.3)}
\]

In spirit of the idea in \cite{HT1}, we establish the following Lemma.

\textbf{Lemma $3.3$}\ \ Let $\Phi_{B}$ be a fixed Schwartz function defined as above and write $u_{B}= \Phi_{B}u$. Then the following estimates hold with constants
 depending on $\Phi_{B}$
\[
\renewcommand{\arraystretch}{1.25}
\begin{array}{lll}
\|u_{B}\|_{H^{\frac{1}{2},2}}\lesssim \|u\|_{\dot{X}^{\frac{1}{2}}_{\zeta}},
\end{array}
\eqno{(3.4)}
\]
\[
\renewcommand{\arraystretch}{1.25}
\begin{array}{lll}
\|u_{B}\|_{H^{1,2}}\lesssim s^{\frac{1}{2}} \|u\|_{\dot{X}^{\frac{1}{2}}_{\zeta}}.
\end{array}
\eqno{(3.5)}
\]

\textbf{Proof}\ \ To prove $(3.4)$, in view of $(3.2)$, it suffices to show that
\[
\renewcommand{\arraystretch}{1.25}
\begin{array}{lll}
\|u_{B}\|_{H^{\frac{1}{2},2}}\lesssim \|u_{B}\|_{{X}^{\frac{1}{2}}_{\zeta}}.
\end{array}
\eqno{(3.6)}
\]

Observing the support of $\Phi_{B}$ is compact, we can find a ball $\tilde{B}\supset \text{spt}\Phi_{B}$ and a cutoff function $\varphi_{\tilde{B}}\in C_0^{\infty}(\tilde{B})$
 satisfying $\varphi_{\tilde{B}}=1$ on $\text{spt}\Phi_{B}$. Then we have $u_{B}=\varphi_{\tilde{B}}u_{B}$. By Fourier transform formula we know
\[
\renewcommand{\arraystretch}{1.25}
\begin{array}{lll}
\hat{u_B}=\hat{\varphi}_{\tilde{B}}\ast \hat{u_B},
\end{array}
\]
where $\hat{\varphi}_{\tilde{B}}$ is a schwartz function.

To obtain $(3.6)$,  by Proposition $3.1$ we just need show that there exists a constant $C$ s.t., $\forall\xi\in \mathbb{R}^n$,
\[
\renewcommand{\arraystretch}{1.25}
\begin{array}{lll}
\displaystyle\int|\hat{\varphi}_{\tilde{B}}(\xi-\eta)|\frac{|\xi|}{|P_{\zeta}(\eta)|+s}d\eta\leq C.
\end{array}
\eqno{(3.7)}
\]

Noting $|\xi|\leq |\xi-\eta|+|\eta|$ and $\hat{\varphi}_{\tilde{B}}$ is a Schwartz function, we only need prove,
\[
\renewcommand{\arraystretch}{1.25}
\begin{array}{lll}
\displaystyle\int|\hat{\varphi}_{\tilde{B}}(\xi-\eta)|\frac{|\eta|}{|P_{\zeta}(\eta)|+s}d\eta\leq C.
\end{array}
\eqno{(3.8)}
\]

On the other hand, if $|\eta|\gg s$ we know $|P_{\zeta}(\eta)|+s\gtrsim |\eta|^2$, which implies
\[
\renewcommand{\arraystretch}{1.25}
\begin{array}{lll}
\dfrac{|\eta|}{|P_{\zeta}(\eta)|+s}\leq C,
\end{array}
\eqno{(3.9)}
\]
where the constant $C$ is independent of $s$. Clearly, $(3.9)$ implies $(3.8)$.

To prove $(3.5)$, observe that, $\forall\eta\in\mathbb{R}^n$,
\[
\renewcommand{\arraystretch}{1.25}
\begin{array}{lll}
\dfrac{|\eta|^2}{|P_{\zeta}(\eta)|+s}\lesssim s.
\end{array}
\eqno{(3.10)}
\]

We can use the same way of the proof of $(3.4)$ to show that
\[
\renewcommand{\arraystretch}{1.25}
\begin{array}{lll}
\|u_{B}\|_{H^{1,2}}\lesssim s^{\frac{1}{2}} \|u\|_{\dot{X}^{\frac{1}{2}}_{\zeta}}.
\end{array}
\eqno{(3.11)}
\]

Then Lemma $3.3$ follows.

Clearly, let $\gamma\in C^1(\bar{\Omega})$ and let $u$ be a solution to $(1.1)$. Then $u$ satisfies
\[
\renewcommand{\arraystretch}{1.25}
\begin{array}{lll}
(-\triangle-A\cdot\nabla)u=0\ \ \ \ \text{in}\  \Omega,
\end{array}
\eqno{(3.12)}
\]
where
\[
\renewcommand{\arraystretch}{1.25}
\begin{array}{lll}
A=\nabla\text{log}\gamma.
\end{array}
\]

The following analysis is from \cite{K1}. For completeness, we will write it down in details. Indeed, we will find the following type of CGO solutions
\[
\renewcommand{\arraystretch}{1.25}
\begin{array}{lll}
u(x,\zeta)=e^{-\frac{\phi_t}{2}}e^{x\cdot\zeta}(1+w(x,\zeta))
\end{array}
\]
of $(3.12)$. Here $\phi_t$ is defined as section $2$ and $\zeta\in C^n\setminus{0}$ satisfies $\zeta\cdot\zeta=0$ which implies that $\text{exp}(x\cdot\zeta)$
is harmonic.

We will decompose
\[
\renewcommand{\arraystretch}{1.25}
\begin{array}{lll}
-\triangle-A\cdot\nabla=-\triangle+(A_t-A)\cdot\nabla-A_t\cdot\nabla
\end{array}
\]
and use the fact that
\[
\renewcommand{\arraystretch}{1.25}
\begin{array}{lll}
(-\triangle-A_t\cdot\nabla)e^{-\frac{\phi_t}{2}}v=e^{-\frac{\phi_t}{2}}(-\triangle+\frac{1}{2}\nabla\cdot A_t+\frac{1}{4}(A_t)^2)v.
\end{array}
\]

Since $A_t=\nabla\phi_t$, it follows that
\[
\renewcommand{\arraystretch}{1.25}
\begin{array}{lll}
(-\triangle-A\cdot\nabla)e^{-\frac{\phi_t}{2}}v&=(-\triangle+(A_t-A)\cdot\nabla-A_t\cdot\nabla)e^{-\frac{\phi_t}{2}}v\\
&=e^{-\frac{\phi_t}{2}}(-\triangle+(A_t-A)\cdot\nabla+q_t)v,
\end{array}
\]
where $q_t=\frac{1}{2}\nabla\cdot A_t-\frac{1}{4}(A_t)^2+\frac{1}{2}A\cdot A_t$.

Hence the equation for $w(x,\zeta)$ is
\[
\renewcommand{\arraystretch}{1.25}
\begin{array}{lll}
(-\triangle_\zeta+(A_t-A)\cdot\nabla_\zeta+q_t)w=(A-A_t)\cdot\zeta-q_t,
\end{array}
\eqno{(3.13)}
\]
where
\[
\renewcommand{\arraystretch}{1.25}
\begin{array}{lll}
-\triangle_\zeta=\triangle+2\zeta\cdot\nabla,
\end{array}
\]
\[
\renewcommand{\arraystretch}{1.25}
\begin{array}{lll}
\nabla_\zeta=\nabla+\zeta.
\end{array}
\]

For given $k\in \mathbb{R}^n$, we set $\zeta_1=s\eta_1+i(\frac{k}{2}+r\eta_2)$, $\zeta_2=-s\eta_1+i(\frac{k}{2}-r\eta_2)$, where $\eta_1, \eta_2\in S^{n-1}$ satisfy
$k\cdot \eta_1=k\cdot \eta_2=\eta_1\cdot \eta_2=0$ and $\frac{|k|^2}{4}+r^2=s^2$. The vectors $\zeta_i$ are chosen so that $\zeta_i\cdot\zeta_i=0$ and $\zeta_1+\zeta_2=k$.

Our goal is to find a sequences $s_{n}$, $\zeta^{(n)}_{i}$ such that $s_{n}\rightarrow\infty$ and $\|w^{(n)}\|_{\dot{X}_{\zeta^{(n)}_{i}}^{\frac{1}{2}}}\rightarrow 0$,
which are the solutions of $(3.13)$ with $t=s_{n}$. In fact we have the following Lemma.

\textbf{Lemma $3.4$}\ \ Let $\gamma$ be $C^{1}(\mathbb{R}^n)$ function with $\gamma>C>0$, $r=1$ outside a ball.
 Then for fixed $k$,  there exists a sequence $\zeta^{(n)}_{i}$ with
$s_{n}\rightarrow\infty$ such that
\[
\renewcommand{\arraystretch}{1.25}
\begin{array}{lll}
\|w^{(n)}\|_{\dot{X}_{\zeta^{(n)}_{i}}^{\frac{1}{2}}}\rightarrow 0,\ \ \ \ \text{as}\ \  s_{n} \rightarrow\infty,
\end{array}
\]
where $w^{(n)}$ is a solution of $(3.13)$ with $t=s_{n}$.

\textbf{Proof}\ \ We want to find the solution $w$ of $(3.13)$. Indeed, we just need find a fixed point of the following operator
\[
\renewcommand{\arraystretch}{1.25}
\begin{array}{lll}
w=\triangle_\zeta^{-1}((A_s-A)\cdot \nabla_\zeta w+q_sw)+\triangle_\zeta^{-1}((A_s-A)\cdot \zeta+q_s)
\end{array}
\]
in a suitable function space, where $\triangle_\zeta^{-1}$ means ${(\triangle_\zeta^{-1}f)}^{\hat{}}=\frac{1}{|P_{\zeta}(\xi)|}\hat{f}(\xi)$.

First of all, we will show that there exists a sequence  $\zeta^{(n)}_{i}$ with $s_{n}\rightarrow\infty$ such that
\[
\renewcommand{\arraystretch}{1.25}
\begin{array}{lll}
\|\triangle_{\zeta_i^{(n)}}^{-1}((A_{s_n}-A)\cdot\zeta_i^{(n)}+q_{s_{n}})\|_{\dot{X}_{\zeta_i^{(n)}}^{\frac{1}{2}}}\longrightarrow 0.
\end{array}
\eqno{(3.14)}
\]

In fact, since $q_{s_n}=\frac{1}{2}\nabla\cdot A_{s_n}-\frac{1}{4}(A_{s_n})^2+\frac{1}{2}A\cdot A_{s_n}$ and  $A_{s_n}$, $A$ is bounded with uniformly compact support, by
Proposition $(3.2)$ we have
\[
\renewcommand{\arraystretch}{1.25}
\begin{array}{lll}
\|(A_{s_n})^2\|_{\dot{X}_{\zeta}^{-\frac{1}{2}}}=\|\Phi_B(A_{s_n})^2\|_{\dot{X}_{\zeta}^{-\frac{1}{2}}}\lesssim
\|(A_{s_n})^2\|_{{X}_{\zeta}^{-\frac{1}{2}}}\lesssim\frac{1}{{s_n}^{\frac{1}{2}}},
\end{array}
\eqno{(3.15)}
\]
\[
\renewcommand{\arraystretch}{1.25}
\begin{array}{lll}
\|A\cdot A_{s_n}\|_{\dot{X}_{\zeta}^{-\frac{1}{2}}}=\|\Phi_B (A\cdot A_{s_n})\|_{\dot{X}_{\zeta}^{-\frac{1}{2}}}\lesssim
\|A\cdot A_{s_n}\|_{{X}_{\zeta}^{-\frac{1}{2}}}\lesssim\frac{1}{{s_n}^{\frac{1}{2}}}.
\end{array}
\eqno{(3.16)}
\]

Now, it's only left to estimate $\|\nabla\cdot A_{s_n}\|_{\dot{X}_{\zeta}^{-\frac{1}{2}}}$. In view of Lemma $3.1$ in \cite{HT1}, we know
\[
\renewcommand{\arraystretch}{1.25}
\begin{array}{lll}
 \displaystyle\int_{S^{n-1}}\frac{1}{s}\displaystyle\int_{s}^{2s}\|\nabla A\|_{\dot{X}_{\zeta_i}^{-\frac{1}{2}}} d\lambda d\eta_1  \longrightarrow 0, \ \ \text{as}\
s\rightarrow\infty.
\end{array}
\]

Then there exist a sequence $\zeta^{(n)}_{i}$ such that
\[
\renewcommand{\arraystretch}{1.25}
\begin{array}{lll}
\|\nabla\cdot A\|_{\dot{X}_{\zeta_i^{(n)}}^{-\frac{1}{2}}}  \longrightarrow 0, \ \ \text{as}\
s_{n}\rightarrow\infty.
\end{array}
\]

Observing $|(\nabla\cdot A_{s_n})^{\hat{}}(\xi)|=|\xi\cdot \hat{A_{s_n}}|=|\xi\cdot\hat{\Psi}(\frac{\xi}{{s_n}})\hat{A}(\xi)|\leq
\|\hat{\Psi}(\frac{\xi}{{s_n}})\|_{L^\infty}|\xi\cdot \hat{A}(\xi)|\lesssim |(\nabla\cdot A)^{\hat{}}(\xi)|$,
then we deduce
\[
\renewcommand{\arraystretch}{1.25}
\begin{array}{lll}
\|\nabla\cdot A_{s_{n}}\|_{\dot{X}_{\zeta_i^{(n)}}^{-\frac{1}{2}}}  \longrightarrow 0, \ \ \text{as}\
s_{n}\rightarrow\infty.
\end{array}
\eqno{(3.17)}
\]

Combining the estimates $(3.15)$, $(3.16)$ and $(3.17)$ we have
\[
\renewcommand{\arraystretch}{1.25}
\begin{array}{lll}
\|q_{s_{n}}\|_{\dot{X}_{\zeta_i^{(n)}}^{-\frac{1}{2}}}  \longrightarrow 0, \ \ \text{as}\
s_{n}\rightarrow\infty.
\end{array}
\]

Secondly, we need show $\|A-A_{s_{n}}\|_{\dot{X}_{\zeta_i^{(n)}}^{-\frac{1}{2}}}=o(\frac{1}{s_{n}})$.

Indeed,  by Lemma $3.1$  again in \cite{HT1} we can find the same sequence  $\zeta^{(n)}_{i}$ with $s_{n}\rightarrow\infty$ such that
\[
\renewcommand{\arraystretch}{1.25}
\begin{array}{lll}
\|A-A_{s_{n}}\|_{\dot{X}_{\zeta_i^{(n)}}^{-\frac{1}{2}}}=\|\Phi_{B}\nabla(\phi-\phi_{s_{n}})\|_{\dot{X}_{\zeta_i^{(n)}}^{-\frac{1}{2}}}\lesssim
\frac{1}{s_n}\|\phi-\phi_{s_{n}}\|_{H^{1, 2}},
\end{array}
\]

From the assumptions of $\gamma$, it follows that $\phi \in C_{c}^{1}(\mathbb{R}^n)$, $\phi_{s_{n}}\in C_{c}^{1}(\mathbb{R}^n)$ with uniformly bounded 
compact support. Hence Lemma $2.1$ implies

\[
\renewcommand{\arraystretch}{1.25}
\begin{array}{lll}
\|A-A_{s_{n}}\|_{L^2}=o(\frac{1}{s_{n}}).
\end{array}
\]

Therefore,  $(3.14)$ holds.

In the following, we will show that  $\triangle_{\zeta_i^{(n)}}^{-1}((A_{s_{n}}-A)\cdot\nabla_{\zeta_i^{(n)}}w+q_{s_{n}}w)$ is a
contract map in $\dot{X}_{\zeta_i^{(n)}}^{-\frac{1}{2}}(\mathbb{R}^n)$ when $s_n$ is large enough. In fact, by corollary $2.1$ in \cite{HT1} we know
\[
\renewcommand{\arraystretch}{1.25}
\begin{array}{lll}
\|(A_{s_{n}}-A)\cdot \zeta_i^{(n)}w\|_{\dot{X}_{\zeta_i^{(n)}}^{\frac{1}{2}}\rightarrow \dot{X}_{\zeta_i^{(n)}}^{-\frac{1}{2}} }\lesssim \|A-A_{s_n}\|_{L_\infty},
\end{array}
\]
\[
\renewcommand{\arraystretch}{1.25}
\begin{array}{lll}
\|q_{s_n}w\|_{\dot{X}_{\zeta_i^{(n)}}^{\frac{1}{2}}\rightarrow \dot{X}_{\zeta_i^{(n)}}^{-\frac{1}{2}} }\lesssim \frac{1}{s_n}\|q_{s_n}\|_{L^\infty}\lesssim
\frac{1}{s_n}(\|A\|_{L^\infty}^{2}+\|D^2\phi_{s_n}\|_{L^\infty}).
\end{array}
\]

On the other hand, by Lemma $2.3$ in \cite{HT1} we know
\[
\renewcommand{\arraystretch}{1.25}
\begin{array}{lll}
\|(A_{s_n}-A)\cdot \nabla w\|_{\dot{X}_{\zeta_i^{(n)}}^{\frac{1}{2}}\rightarrow \dot{X}_{\zeta_i^{(n)}}^{-\frac{1}{2}} }\lesssim \|A-A_{s_n}\|_{L^\infty}.
\end{array}
\]

Then, by Lemma $2.1$ it follows that
\[
\renewcommand{\arraystretch}{1.25}
\begin{array}{lll}
\|\triangle_{\zeta_i^{(n)}}^{-1}((A_{s_{n}}-A)\cdot\nabla_{\zeta_i^{(n)}}w+q_{s_{n}}w)\|_{\dot{X}_{\zeta_i^{(n)}}^{\frac{1}{2}}\rightarrow \dot{X}_{\zeta_i^{(n)}}^{\frac{1}{2}}}
< 1,
\end{array}
\]
when $s_n$ is large enough.

Using contract mapping theorem we know there exist a sequence $\zeta_i^{(n)}$ such that
\[
\renewcommand{\arraystretch}{1.25}
\begin{array}{lll}
\|w^{(n)}\|_{\dot{X}_{\zeta^{(n)}_{i}}^{\frac{1}{2}}}\rightarrow 0,\ \ \ \ \text{as}\ \  s_{n} \rightarrow\infty,
\end{array}
\]
where $w^{(n)}$ is a solution of $(3.13)$ with $t=s_{n}$.

\textbf{Lemma $3.5$}\ \ Let $w^{(n)}$ be chosen as in Lemma $3.4$, then the following estimates hold
\[
\renewcommand{\arraystretch}{1.25}
\begin{array}{lll}
\|w^{(n)}\|_{L^{2}(\Omega)}\lesssim \frac{1}{s_{n}^{\frac{1}{2}}},
\end{array}
\]
\[
\renewcommand{\arraystretch}{1.25}
\begin{array}{lll}
\|w^{(n)}\|_{H^{1,2}(\Omega)}\lesssim s_{n}^{\frac{1}{2}},
\end{array}
\]
\[
\renewcommand{\arraystretch}{1.25}
\begin{array}{lll}
\|w^{(n)}\|_{H^{2,2}(\Omega)}\lesssim s_{n}^{\frac{3}{2}}.
\end{array}
\]

\textbf{Proof}\ \ From the construction of $w^{(n)}$, by Proposition $3.2$ and Lemma $3.3$  we know
\[
\renewcommand{\arraystretch}{1.25}
\begin{array}{lll}
\|w^{(n)}\|_{L^{2}(\Omega)}\leq \|\Phi_{B}w^{(n)}\|_{L^{2}}\lesssim \frac{1}{s_{n}^{\frac{1}{2}}}\|w^{(n)}\|_{\dot{X}_{\zeta^{(n)}_{i}}^{\frac{1}{2}}}
\end{array}
\eqno{(3.18)}
\]
and
\[
\renewcommand{\arraystretch}{1.25}
\begin{array}{lll}
\|\Phi_{B}w^{(n)}\|_{H^{1,2}}\lesssim s_{n}^{\frac{1}{2}}\|w^{(n)}\|_{\dot{X}_{\zeta^{(n)}_{i}}^{\frac{1}{2}}}\lesssim s_{n}^{\frac{1}{2}}.
\end{array}
\eqno{(3.19)}
\]

In view of the definition of $\Phi_{B}$, we know, for $\Omega\subset B_0 \subset B$, $(3.19)$ implies
\[
\renewcommand{\arraystretch}{1.25}
\begin{array}{lll}
\|\Phi_{B}w^{(n)}\|_{H^{1,2}(B)}\lesssim s_{n}^{\frac{1}{2}}.
\end{array}
\eqno{(3.20)}
\]

Thus $w^{(n)}$ is a weak solution of the following equation
\[
\renewcommand{\arraystretch}{1.25}
\begin{array}{lll}
-\triangle w^{(n)}+(2\zeta^{(n)}_{i}+(A_{s_n}-A))\cdot\nabla w^{(n)}&+(\zeta^{(n)}_{i}\cdot(A_{s_n}-A)+q_{s_n})w^{(n)}\\
&=(A-A_{s_n})\cdot \zeta^{(n)}_{i}-q_{s_n},\ \ \ \ \ \ \ \ \ \ \ \ \ \ \ \text{in}\ \  B.
\end{array}
\]

By the classical interior estimate for Laplace equation we deduce
 \[
\renewcommand{\arraystretch}{1.25}
\begin{array}{lll}
\|w^{(n)}\|_{H^{2,2}(\Omega)}&\lesssim \|w^{(n)}\|_{H^{1,2}(B_0)} +\|(2\zeta^{(n)}_{i}+(A_{s_n}-A))\cdot\nabla w^{(n)}\|_{L^2(B_0)}\\
&+\|(\zeta^{(n)}_{i}\cdot(A_{s_n}-A)+q_{s_n})w^{(n)}\|_{L^2(B_0)}+\|(A-A_{s_n})\cdot \zeta^{(n)}_{i}-q_{s_n}\|_{L^2(B_0)}\\
&\lesssim \|w^{(n)}\|_{H^{1,2}(B_0)}+s_n\|w^{(n)}\|_{H^{1,2}(B_0)}+s_n\|w^{(n)}\|_{L^{2}(B_0)}+s_n,
\end{array}
\eqno{(3.21)}
\]
where we use the relations $\|A_{s_n}\|_{L^\infty}\leq \|A\|_{L^\infty}$ and $\|D^2\Phi_{s_n}\|_{L^\infty}=o(s_n)$.

Combining $(3.20)$ and $(3.21)$ we have
\[
\renewcommand{\arraystretch}{1.25}
\begin{array}{lll}
\|w^{(n)}\|_{H^{2,2}(\Omega)}\lesssim s_{n}^{\frac{3}{2}}.
\end{array}
\]

The proof is complete.

\bigbreak
\begin{center}
\textbf{\Large {\S 4. Carleman Estimate for CGO Solutions }}
\end{center}

In this section, we will introduce Carleman estimate for CGO solutions. The first result is taken from \cite{K1}.

\textbf{Proposition $4.1$}\ \ Let $\xi\in S^{n-1}$ and  suppose $u\in H^2(\Omega)$. Then there exists a constant $t_0>0$ such that
for $t\geq t_0$, we have the estimate
\[
\renewcommand{\arraystretch}{1.25}
\begin{array}{lll}
C(t^2\|u\|_{L^2(\Omega)}^{2}+\|\nabla u\|_{L^2(\Omega)}^2)-C't^2\displaystyle\int_{\partial \Omega}|u|^2ds-C''\displaystyle\int_{\partial \Omega}\bar{u}\partial_{\nu}uds+\\
\displaystyle\int_{\partial \Omega}4t\text{Re}(\partial_{\nu}u\partial_{\eta}\bar{u})-2t(\nu\cdot\eta)|\nabla u|^2+2t^3(\nu\cdot\eta)|u|^2ds\leq
\|e^{-x\cdot t\eta}(-\triangle)(e^{x\cdot t\eta}u)\|_{L^2(\Omega)}^2.
\end{array}
\eqno{(4.1)}
\]

In this paper, we need the following type of Carleman estimate.

\textbf{Lemma $4.2$}\ \ Let $\eta\in S^{n-1}$ and suppose $u\in H^2(\Omega)$. Then there exists a constant $\delta$ such that for $\gamma\in C^1(\bar{\Omega})$,  we have
 \[
 \begin{array}{lll}
C(t^2\|u\|_{L^2(\Omega)}^{2}+\|\nabla u\|_{L^2(\Omega)}^2)-C't^2\displaystyle\int_{\partial \Omega}|u|^2ds-C''\displaystyle\int_{\partial \Omega}\bar{u}\partial_{\nu}uds\\
+\displaystyle\int_{\partial \Omega}4t\text{Re}(\partial_{\nu}u\partial_{\eta}\bar{u})-2t(\nu\cdot\eta)|\nabla u|^2+2t^3(\nu\cdot\eta)|u|^2ds\\
\leq
\|e^{-x\cdot t\eta}(-\triangle+(A_t-A)\cdot\nabla+q_t)(e^{x\cdot t\eta}u)\|_{L^2(\Omega)}^2
\end{array}
\eqno{(4.2)}
\]
for  $t\geq t{(\gamma)}>0$, where $t{(\gamma)}$ is a constant only depending on $\gamma$.

\textbf{Proof}\ \ Since
 \[
 \begin{array}{lll}
\|e^{-x\cdot t\eta}(-\triangle)(e^{x\cdot t\eta}u)\|_{L^2(\Omega)}^2&\leq
2\|e^{-x\cdot t\eta}(-\triangle+(A_t-A)\cdot\nabla+q_t)(e^{x\cdot t\eta}u)\|_{L^2(\Omega)}^2\\
&+\|e^{-x\cdot t\eta}((A_t-A)\cdot\nabla+q_t)(e^{x\cdot t\eta}u)\|_{L^2(\Omega)}^2\\
&\leq 2\|e^{-x\cdot t\eta}(-\triangle+(A_t-A)\cdot\nabla+q_t)(e^{x\cdot t\eta}u)\|_{L^2(\Omega)}^2\\
&+2\|A-A_t\|_{L^\infty}^2\|\nabla u\|_{L^2(\Omega)}^2
+ 2t^2\|A-A_t\|_{L^\infty}^2\|u\|_{L^2(\Omega)}^2\\
&+2\|q_t\|_{L^\infty}^2\|u\|_{L^2(\Omega)}^2.
\end{array}
\eqno{(4.3)}
\]

Since $\gamma\in C^1(\bar{\Omega})$, from Lemma $2.1$, we know
$\|A-A_t\|_{L^\infty}=o(1)$ and $\|q_t\|_{L^\infty}=o(t)$. Hence, if $t$ is large enough,  the superfluous terms on
the right hand in $(4.3)$ can be absorbed by the first term on the left hand in $(4.1)$. Then $(4.2)$ holds.

\bigbreak
\begin{center}
\textbf{\Large {\S 5. The Uniqueness Proof }}
\end{center}

First, we introduce a boundary integral identity which is from \cite{K1}.

\textbf{Proposition $5.1$}\ \ Suppose $\gamma_{i}\in C^1(\bar{\Omega}), i=1, 2.$ and $u_1, u_2\in H^1(\Omega)$
satisfy $\nabla\cdot(\gamma_i\nabla u_i)=0$ in $\Omega$. Suppose further that $\tilde{ u}_1\in H^1(\Omega)$
satisfies $\nabla\cdot(\gamma_1\nabla \tilde{ u}_1)=0$ with $\tilde{ u}_1=u_2$ on $\partial \Omega$. Then
 \[
 \begin{array}{lll}
\displaystyle\int_{\Omega}(\gamma_1^{\frac{1}{2}}\nabla \gamma_2^{\frac{1}{2}}-\gamma_2^{\frac{1}{2}}\nabla \gamma_1^{\frac{1}{2}})\cdot \nabla (u_1 u_2)dx=
\displaystyle\int_{\partial \Omega}\gamma_1 \partial_{\nu}(\tilde{ u}_1-u_2)u_1ds,
\end{array}
\]
where the integral on the boundary is understood in the sense of the dual pairing between $H^{\frac{1}{2}}(\partial \Omega)$ and $H^{-\frac{1}{2}}(\partial \Omega)$.

\textbf{Remark $5.2$}\ \ Proposition $5.1$ also holds for $\gamma_i\in W^{1,\infty}(\bar{\Omega})$.

we will make use of boundary integral identity to show the  uniqueness of the  Calder\'{o}n problem with partial data. For this goal, fix $k\in \mathbb{R}^n $ with
 $k\cdot\eta=0$ and choose $l^{(n)}\in \mathbb{R}^n $ with $l^{(n)}\cdot \eta =l^{(n)}\cdot k=0$ and $\frac{|k|^2}{4}+|l^{(n)}|^2=s_n^2.$ Further define 
$\zeta_1^{(n)}=s_n\eta
+i(\frac{k}{2}+l^{(n)})$ and $\zeta_2^{(n)}=-s_n\eta
+i(\frac{k}{2}-l^{(n)})$,   and
let $\phi_i=\nabla \text{log}\gamma_i$ and $u^{(n)}_i=e^{-\frac{\phi_{is_n}}{2}}e^{x\cdot\zeta_i^{(n)}}(1+w_i^{(n)})$ be
CGO solutions to $\nabla\cdot(\gamma_i\nabla u_i)=0$, where $w_i^{(n)}$ is chosen as Lemma $3.4$. Finally take $\tilde {u}_1^{(n)}$ as the solution
 to $\nabla\cdot(\gamma_1 \nabla \tilde {u}_1^{(n)})$ with $\tilde {u}_1^{(n)}=u_2^{(n)}$ on $\partial \Omega$. Hence Proposition $5.1$ implies
 \[
 \begin{array}{lll}
\displaystyle\int_{\Omega}(\gamma_1^{\frac{1}{2}}\nabla \gamma_2^{\frac{1}{2}}-\gamma_2^{\frac{1}{2}}\nabla \gamma_1^{\frac{1}{2}})\cdot \nabla (u_1^{(n)} u_2^{(n)})dx=
\displaystyle\int_{\partial \Omega}\gamma_1 \partial_{\nu}(\tilde{u}_1^{(n)}-u_2^{(n)})u_1^{(n)}ds.
\end{array}
\eqno{(5.1)}
\]

We will first prove
 \[
 \begin{array}{lll}
\lim\limits_{n\rightarrow\infty}\displaystyle\int_{\partial \Omega}\gamma_1 \partial_{\nu}(\tilde{u}_1^{(n)}-u_2^{(n)})u_1^{(n)}ds=0.
\end{array}
\eqno{(5.2)}
\]

 The proof of $(5.2)$ is divided into  the following three Lemmas. For simplicity, we will not write the superscripts of
$u_i^{(n)}$,$\tilde {u}_1^{(n)}$, $w_i^{(n)}$ and $\zeta_i^{(n)}$ and the subscript of $s_n$ again unless otherwise particularly specified.

Introduce the function
\[
 \begin{array}{lll}
u=e^{\frac{\phi_{1s}}{2}}\tilde{u}_1-e^{\frac{\phi_{2s}}{2}}{u}_2=u_0+\delta u,
\end{array}
\]
where
\[
 \begin{array}{lll}
u_0=e^{\frac{\phi_{1s}}{2}}(\tilde{u}_1-u_2),
\end{array}
\]
\[
 \begin{array}{lll}
\delta u=(e^{\frac{\phi_{1s}}{2}}-e^{\frac{\phi_{2s}}{2}})u_2.
\end{array}
\]

\textbf{Lemma $5.3$}\ \ Suppose the assumptions of Theorem $1.1$ hold and define $u_2$ as above. Then
\[
 \begin{array}{lll}
\displaystyle\int_{\partial \Omega}e^{-x\cdot 2s\eta}|\nabla \delta u|^2ds=o(1),
\end{array}
\]
\[
 \begin{array}{lll}
\displaystyle\int_{\partial \Omega}e^{-x\cdot 2s\eta}| \delta u|^2ds=o(\frac{1}{s^2}),
\end{array}
\]
as $s\rightarrow \infty$.

\textbf{Proof}\ \ From the definition of $\delta u$, we know
\[
 \begin{array}{lll}
\displaystyle\int_{\partial \Omega}e^{-x\cdot 2s\eta}|\nabla \delta u|^2ds &\lesssim \displaystyle\int_{\partial \Omega}e^{-x\cdot 2s\eta}
|\nabla(e^{\frac{\phi_{1s}}{2}}-e^{\frac{\phi_{2s}}{2}})|^2|u_2|^2ds\\
&+\displaystyle\int_{\partial \Omega}e^{-x\cdot 2s\eta}|e^{\frac{\phi_{1s}}{2}}-e^{\frac{\phi_{2s}}{2}}|^2|\nabla u_2|^2ds=:I+II.
\end{array}
\]

For $I$, from the construction of $u_2$, we have
\[
 \begin{array}{lll}
I\lesssim \displaystyle\int_{\partial \Omega}|\nabla(e^{\frac{\phi_{1s}}{2}}-e^{\frac{\phi_{2s}}{2}})|^2|1+w_2|^2ds,
\end{array}
\eqno{(5.3)}
\]

Since $\gamma_1$\textbar$ _{\partial \Omega}=\gamma_2$\textbar$ _{\partial \Omega}$, $\partial_{\nu}\gamma_1$\textbar$ _{\partial \Omega}
=\partial_{\nu}\gamma_2$\textbar$ _{\partial \Omega}$, by mean value theorem  $(5.3)$ implies
\[
 \begin{array}{lll}
I&\lesssim \displaystyle\int_{\partial \Omega}|\nabla(e^{\frac{\phi_{1s}}{2}}-\gamma_1^{\frac{1}{2}} )|^2ds
+\displaystyle\int_{\partial \Omega}|\nabla(e^{\frac{\phi_{2s}}{2}}-\gamma_2^{\frac{1}{2}} )|^2ds\\
&+\displaystyle\int_{\partial \Omega}|w_2|^2ds\|\nabla(e^{\frac{\phi_{1s}}{2}}-\gamma_1^{\frac{1}{2}} )\|_{L^\infty}^2
+\displaystyle\int_{\partial \Omega}|w_2|^2ds\|\nabla(e^{\frac{\phi_{2s}}{2}}-\gamma_2^{\frac{1}{2}} )\|_{L^\infty}^2\\
&\lesssim \displaystyle\int_{\partial \Omega}|\nabla \phi_{1s}-\nabla \phi_1|^2ds+\displaystyle\int_{\partial \Omega}| \phi_{1s}-\phi_1|^2ds\\
&+\displaystyle\int_{\partial \Omega}|\nabla \phi_{2s}-\nabla \phi_2|^2ds+\displaystyle\int_{\partial \Omega}| \phi_{2s}-\phi_2|^2ds\\
&+\displaystyle\int_{\partial \Omega}|w_2|^2ds(\|\nabla(\phi_{1s}-\phi_1)\|_{L^\infty}^2+\|\phi_{1s}-\phi_1\|_{L^\infty}^2)\\
&+\displaystyle\int_{\partial \Omega}|w_2|^2ds(\|\nabla(\phi_{2s}-\phi_2)\|_{L^\infty}^2+\|\phi_{2s}-\phi_2\|_{L^\infty}^2).
\end{array}
\]

Hence it follows that 
\[
 \begin{array}{lll}
I&\lesssim \|A_{1s}-A_1\|_{L^{\infty}}^{2}+\|\phi_{1s}-\phi_1\|_{L^{\infty}}^2+\|A_{2s}-A_2\|_{L^{\infty}}^{2}+\|\phi_{2s}-\phi_2\|_{L^{\infty}}^2\\
&+\displaystyle\int_{\partial \Omega}|w_2|^2ds(\|A_{1s}-A_1\|_{L^{\infty}}^2+\|\phi_{1s}-\phi_1\|_{L^{\infty}}^2)\\
&+\displaystyle\int_{\partial \Omega}|w_2|^2ds(\|A_{2s}-A_2\|_{L^{\infty}}^2+\|\phi_{2s}-\phi_2\|_{L^{\infty}}^2).
\end{array}
\eqno{(5.4)}
\]

On the other hand, by Lemma $2.2$ and Lemma  $3.5$ we have 
\[
 \begin{array}{lll}
\displaystyle\int_{\partial \Omega}|w_2|^2ds \lesssim 1.
\end{array}
\eqno{(5.5)}
\]

By Lemma $2.1$, $(5.4)$ and $(5.5)$ implies
\[
 \begin{array}{lll}
I=o(1),\ \ \ \ \text{as}\ \  s\rightarrow \infty.
\end{array}
\]

For $II$, from the definition of $u_2$, we know
\[
 \begin{array}{lll}
II &\lesssim \displaystyle\int_{\partial \Omega}|e^{\frac{\phi_{1s}}{2}}-\gamma_1^{\frac{1}{2}}|^2|\nabla w_2|^2ds+
\displaystyle\int_{\partial \Omega}|e^{\frac{\phi_{2s}}{2}}-\gamma_2^{\frac{1}{2}}|^2|\nabla w_2|^2ds\\
&+s^2 \displaystyle\int_{\partial \Omega}|e^{\frac{\phi_{1s}}{2}}-\gamma_1^{\frac{1}{2}}|^2|1+ w_2|^2ds+s^2 \displaystyle\int_{\partial \Omega}|e^{\frac{\phi_{2s}}{2}}-\gamma_2^{\frac{1}{2}}|^2|1+ w_2|^2ds.
\end{array}
\]

Hence it follows that 
\[
 \begin{array}{lll}
II &\lesssim \|\phi_{1s}-\phi_1\|_{L^\infty}^2\displaystyle\int_{\partial \Omega}|\nabla w_2|^2ds+\|\phi_{2s}-\phi_2\|_{L^\infty}^2\displaystyle\int_{\partial \Omega}|\nabla w_2|^2ds\\
&+s^2\|\phi_{1s}-\phi_1\|_{L^\infty}^2\displaystyle\int_{\partial \Omega}1+|w_2|^2ds+s^2\|\phi_{2s}-\phi_2\|_{L^\infty}^2\displaystyle\int_{\partial \Omega}1+|w_2|^2ds.
\end{array}
\eqno{(5.6)}
\]

On the other hand, by Lemma $2.2$ and Lemma $3.5$ we have
\[
 \begin{array}{lll}
\displaystyle\int_{\partial \Omega}|\nabla w_2|^2ds \lesssim s^2.
\end{array}
\eqno{(5.7)}
\]

By Lemma $2.1$, $(5.5)$, $(5.6)$ and $(5.7)$ implies
\[
 \begin{array}{lll}
II=o(1),\ \ \ \ \text{as}\ \  s\rightarrow \infty.
\end{array}
\]

Hence, we prove
\[
 \begin{array}{lll}
\lim\limits_{s\rightarrow\infty}\displaystyle\int_{\partial \Omega}e^{-x\cdot 2s\eta}|\nabla \delta u|^2ds=0.
\end{array}
\]

Noting
\[
 \begin{array}{lll}
\displaystyle\int_{\partial \Omega}e^{-x\cdot 2s\eta}| \delta u|^2ds\lesssim \displaystyle\int_{\partial \Omega}
|e^{\frac{\phi_{1s}}{2}}-e^{\frac{\phi_{2s}}{2}}|^2(1+|w_2|^2)dx
\end{array}
\]
and using the above estimates, we can easily obtain
\[
 \begin{array}{lll}
\displaystyle\int_{\partial \Omega}e^{-x\cdot 2s\eta}| \delta u|^2ds=o(\frac{1}{s^2}),\ \ \ \ \text{as}\ s\rightarrow\infty.
\end{array}
\]

The proof is complete.

\textbf{Lemma $5.4$}\ \ Suppose the assumptions of Theorem $1.1$ hold and define $u_1$,$u_2$ and $\tilde{u}_1$ as above. Then
\[
 \begin{array}{lll}
\displaystyle\int_{\partial \Omega_{+,\varepsilon}}e^{-x\cdot 2s\eta}| \partial_{\nu} u|^2ds=o(1),
\end{array}
\]
as $s\rightarrow \infty$.

\textbf{Proof}\ \ Letting $v=e^{-x\cdot s\eta}u$ and recalling Carleman estimate for $v$ we have
 \[
 \begin{array}{lll}
C(s^2\|v\|_{L^2(\Omega)}^{2}+\|\nabla v\|_{L^2(\Omega)}^2)
-C's^2\displaystyle\int_{\partial \Omega}|v|^2ds-C''\displaystyle\int_{\partial \Omega}\bar{v}\partial_{\nu}vds\\
+\displaystyle\int_{\partial \Omega}4s\text{Re}(\partial_{\nu}v\partial_{\eta}\bar{v})-2s(\nu\cdot\eta)|\nabla v|^2+2s^3(\nu\cdot\eta)|v|^2ds\\
\leq
\|e^{-x\cdot s\eta}(-\triangle+(A_{1s}-A_{1})\cdot\nabla+q_{1s})u\|_{L^2(\Omega)}^2.
\end{array}
\eqno{(5.8)}
\]

Hence
 \[
 \begin{array}{lll}
\displaystyle\int_{\partial \Omega}4\text{Re}(\partial_{\nu}v\partial_{\eta}\bar{v})
-2(\nu\cdot\eta)|\nabla v|^2ds\lesssim s\displaystyle\int_{\partial \Omega}|v|^2ds+\frac{1}{s}|\displaystyle\int_{\partial \Omega}\bar{v}\partial_{\nu}vds|\\
+s^2\displaystyle\int_{\partial \Omega}|\nu\cdot\eta||v|^2ds+\frac{1}{s}\|e^{-x\cdot s\eta}(-\triangle+(A_{1s}-A_1)\cdot\nabla+q_{1s})u\|_{L^2(\Omega)}^2\\
=:I+II+III+IV.
\end{array}
\eqno{(5.9)}
\]

Now we first estimate the terms on the left hand of $(5.9)$.
 \[
 \begin{array}{lll}
\text{Left of}\  (5.9)&=\displaystyle\int_{\partial \Omega_{+,\epsilon}}4\text{Re}(\partial_{\nu}v\partial_{\eta}\bar{v})-2(\nu\cdot\eta)|\nabla v|^2ds\\
&+\displaystyle\int_{\partial \Omega_{-,\epsilon}}4\text{Re}(\partial_{\nu}v\partial_{\eta}\bar{v})-2(\nu\cdot\eta)|\nabla v|^2ds=:V+VI.
\end{array}
\]

For $V$, since $\tilde{u}_1$\textbar$_{\partial \Omega}=u_2$\textbar$_{\partial \Omega}$, recalling the definition of $V$ and using Young
inequality we have
 \[
 \begin{array}{lll}
V &\geq \displaystyle\int_{\partial \Omega_{+,\epsilon}}e^{-x\cdot2s\eta}\nu\cdot\eta(\frac{\partial u}{\partial \nu})^2ds
-C\displaystyle\int_{\partial \Omega_{+,\epsilon}}e^{-x\cdot2s\eta}|\nabla\delta u|^2ds\\
&-Cs^2\displaystyle\int_{\partial \Omega_{+,\epsilon}}e^{-x\cdot2s\eta}|\delta u|^2ds.
\end{array}
\eqno{(5.10)}
\]

For $VI$, observing $\tilde{u}_1$\textbar$_{\partial \Omega}=u_2$\textbar$_{\partial \Omega}$
and $\frac{\partial \tilde{u}_1}{\partial \nu}$\textbar$_{\partial \Omega_{-,\epsilon}}=\frac{\partial u_2}{\partial \nu}$\textbar$_{\partial \Omega_{-,\epsilon}}$
we deduce
 \[
 \begin{array}{lll}
VI&\lesssim \displaystyle\int_{\partial \Omega_{-,\epsilon}}|\nabla v|^2ds \lesssim s^2\displaystyle\int_{\partial \Omega_{-,\epsilon}}e^{-x\cdot2s\eta}|u|^2ds
+\displaystyle\int_{\partial \Omega_{-,\epsilon}}e^{-x\cdot2s\eta}|\nabla u|^2ds\\
&\lesssim s^2\displaystyle\int_{\partial \Omega_{-,\epsilon}}e^{-x\cdot2s\eta}|\delta u|^2ds
+\displaystyle\int_{\partial \Omega_{-,\epsilon}}e^{-x\cdot2s\eta}|\nabla\delta u|^2ds.
\end{array}
\eqno{(5.11)}
\]

Combing $(5.10)$ and $(5.11)$ we deduce
 \[
 \begin{array}{lll}
\text{Left of $(5.9)$}&\geq \displaystyle\int_{\partial \Omega_{+,\epsilon}}e^{-x\cdot2s\eta}\nu\cdot\eta(\frac{\partial u}{\partial \nu})^2ds\\
&-Cs^2\displaystyle\int_{\partial \Omega }e^{-x\cdot2s\eta}|\delta u|^2ds
-C\displaystyle\int_{\partial \Omega } e^{-x\cdot 2s\eta}|\nabla\delta u|^2ds\\
&\geq \epsilon\displaystyle\int_{\partial \Omega_{+,\epsilon}}e^{-x\cdot2s\eta}(\frac{\partial u}{\partial \nu})^2ds-o(1),
\end{array}
\eqno{(5.12)}
\]
as $s\rightarrow \infty$, where we use Lemma $5.3$ in the second inequality.

Now we will deal with the terms of on the right hand of $(5.9)$. For $I$ and $III$, by Lemma $5.3$ we have
 \[
 \begin{array}{lll}
I\lesssim s \displaystyle\int_{\partial \Omega } e^{-x\cdot 2s\eta}|u|^2ds\lesssim s \displaystyle\int_{\partial \Omega } e^{-x\cdot 2s\eta}|\delta u|^2ds=o(\frac{1}{s})
\end{array}
\]
and
 \[
 \begin{array}{lll}
III\lesssim s^2 \displaystyle\int_{\partial \Omega } e^{-x\cdot 2s\eta}|u|^2ds\lesssim s^2 \displaystyle\int_{\partial \Omega } e^{-x\cdot 2s\eta}|\delta u|^2ds=o(1).
\end{array}
\]

From the definition of $v$,  we know
\[
 \begin{array}{lll}
\frac{1}{s}\displaystyle\int_{\partial \Omega }\bar{v}\partial_{\nu}v ds&=
\frac{1}{s}\displaystyle\int_{\partial \Omega }e^{-x\cdot s\eta}\bar{\delta u}\partial_{\nu}(e^{-x\cdot s\eta}u)ds\\
&=\frac{1}{s}\displaystyle\int_{\partial \Omega }e^{-x\cdot s\eta}\bar{\delta u}e^{-x\cdot s\eta}-s\nu\cdot\eta\delta u ds
+\frac{1}{s}\displaystyle\int_{\partial \Omega }e^{-x\cdot s\eta}\bar{\delta u}e^{-x\cdot s\eta}\partial_{\nu}u ds.
\end{array}
\]

Observing $\tilde{u}_1$\textbar$_{\partial \Omega}=u_2$\textbar$_{\partial \Omega}$
and $\frac{\partial \tilde{u}_1}{\partial \nu}$\textbar$_{\partial \Omega_{-,\epsilon}}=\frac{\partial u_2}{\partial \nu}$\textbar$_{\partial \Omega_{-,\epsilon}}$ and
 using Young inequality we have, for any $\epsilon_0>0$,
\[
 \begin{array}{lll}
II&\leq \epsilon_0 \displaystyle\int_{\partial \Omega_{+,\epsilon}}e^{-x\cdot2s\eta}(\frac{\partial u}{\partial \nu})^2ds
+C\displaystyle\int_{\partial \Omega }e^{-x\cdot2s\eta}|\delta u|^2ds+C\frac{1}{s^2}\displaystyle\int_{\partial \Omega}e^{-x\cdot2s\eta}|\nabla\delta u|^2ds\\
&\leq \epsilon_0 \displaystyle\int_{\partial \Omega_{+,\epsilon}}e^{-x\cdot2s\eta}(\frac{\partial u}{\partial \nu})^2ds+o(\frac{1}{s^2}),
\end{array}
\]
as $s\rightarrow \infty$, where we use Lemma $5.3$ in the second inequality.

Finally, we estimate $IV$.
\[
 \begin{array}{lll}
IV&\leq \frac{1}{s}\displaystyle\int_{\Omega}e^{-x\cdot 2s\eta}|(-\triangle+(A_{1s}-A_{1})\cdot\nabla+q_{1s})(e^{\frac{\phi_{2s}}{2}}u_2)|^2dx\\
&\leq \frac{1}{s}\displaystyle\int_{\Omega}e^{-x\cdot 2s\eta}|((A_{1s}-A_{1})-(A_{2s}-A_{2}))\cdot \nabla(e^{x\cdot\zeta_2}(1+w_2))
+(q_{1s}-q_{2s})e^{x\cdot\zeta_2}(1+w_2)|^2dx,
\end{array}
\eqno{(5.13)}
\]
where we use the equations
\[
 \begin{array}{lll}
(-\triangle+(A_{is}-A_{i})\cdot \nabla+q_{is})(e^{-\frac{\phi_{is}}{2}}u_i)=0.
\end{array}
\]

$(5.13)$ implies
\[
 \begin{array}{lll}
IV&\lesssim s\displaystyle\int_{\Omega}|A_2-A_{2s}|^2+|A_1-A_{1s}|^2dx+s\displaystyle\int_{\Omega}(|A_2-A_{2s}|^2+|A_1-A_{1s}|^2)|w_2|^2dx\\
&+\frac{1}{s}\displaystyle\int_{\Omega}(|A_2-A_{2s}|^2+|A_1-A_{1s}|^2)|\nabla w_2|^2dx
+\frac{1}{s}\displaystyle\int_{\Omega}|q_{2s}-q_{1s}|^2dx\\
&+\frac{1}{s}\displaystyle\int_{\Omega}|q_{2s}-q_{1s}|^2|w_2|^2dx
=:a+b+c+d+e.
\end{array}
\]

By Lemma $2.1$  we know
\[
 \begin{array}{lll}
a=o(1),\ \ \ \ d=o(1),\ \ \ \ \text{as}\ s\rightarrow\infty.
\end{array}
\]

For $b$, by proposition $3.2$ and Lemma $3.4$ we have
\[
 \begin{array}{lll}
b&\leq (\|A_2-A_{2s}\|_{L^\infty}^2+\|A_1-A_{1s}\|_{L^\infty}^2)s\displaystyle\int_{\Omega}|w_2|^
2dx\\
&\lesssim s \displaystyle\int_{\Omega}|\Phi_B w_2|^2dx\\
&\lesssim \|w_2\|_{\dot{X}_{\zeta_2}^{\frac{1}{2}}}\longrightarrow 0\ \ \ \ \text{as}\ s\rightarrow \infty.
\end{array}
\]

Observing $\|q_{2s}-q_{1s}\|_{L^\infty}=o(s)$ and using the same way for $b$ we can obtain
\[
 \begin{array}{lll}
e=o(1),\ \ \ \ \text{as}\ s\rightarrow\infty.
\end{array}
\]

Finally, for $c$, using Lemma $2.1$ and Lemma $3.3$ we have
\[
 \begin{array}{lll}
c\lesssim \frac{1}{s}\displaystyle\int_{\Omega}|\nabla w_2|^2 dx \lesssim \frac{1}{s}\displaystyle\int_{\Omega}|\nabla(\Phi_B w_2)|^2 dx\lesssim
 \|w_2\|_{\dot{X}_{\zeta_2}^{\frac{1}{2}}}\longrightarrow 0,\ \ \text{as}\ s\rightarrow\infty.
\end{array}
\]

From the above analysis we obtain
\[
 \begin{array}{lll}
IV=o(1),\ \ \ \ \text{as}\ s\rightarrow\infty.
\end{array}
\]

Combining $(5.9)$, $(5.12)$ and the estimates of $I$, $II$, $III$ and $IV$ we deduce
\[
 \begin{array}{lll}
\lim\limits_{s\rightarrow \infty}\displaystyle\int_{\partial \Omega_{+,\epsilon}}e^{-x\cdot 2s\eta}| \partial_{\nu} u|^2ds=0.
\end{array}
\]

\textbf{Lemma $5.5$}\ \ Suppose the assumptions of Theorem $1.1$ hold and define $u_1$,$u_2$ and $\tilde{u}_1$ as above. Then
\[
 \begin{array}{lll}
\lim\limits_{s\rightarrow\infty}\displaystyle\int_{\partial \Omega}\gamma_1 \partial_{\nu}(\tilde{u}_1-u_2)u_1ds=0.
\end{array}
\]

\textbf{Proof}\ \ From the fact  $\frac{\partial \tilde{u}_1}{\partial \nu}$\textbar$_{\partial \Omega_{-,\epsilon}}
=\frac{\partial u_2}{\partial \nu}$\textbar$_{\partial \Omega_{-,\epsilon}}$ and $(5.5)$ Cauchy-Schwarz inequality gives 
\[
 \begin{array}{lll}
|\displaystyle\int_{\partial \Omega}\gamma_1 \partial_{\nu}(\tilde{u}_1-u_2)u_1ds|&\lesssim \displaystyle\int_{\partial \Omega_{+,\epsilon}}
e^{-x\cdot 2s\eta}|\partial_{\nu}(\tilde{u}_1-u_2)|^2ds \|1+w_1\|_{L^2(\partial \Omega)}^2.\\
&\lesssim \displaystyle\int_{\partial \Omega_{+,\epsilon}}
e^{-x\cdot 2s\eta}|\partial_{\nu}(\tilde{u}_1-u_2)|^2ds.
\end{array}
\]

To estimate $\int_{\partial \Omega_{+,\epsilon}}
e^{-x\cdot 2s\eta}|\partial_{\nu}(\tilde{u}_1-u_2)|^2ds$, from the definition of $u$, $u_0$ and $\delta u$ we deduce
\[
 \begin{array}{lll}
\displaystyle\int_{\partial \Omega_{+,\epsilon}}
e^{-x\cdot 2s\eta}|\partial_{\nu}(\tilde{u}_1-u_2)|^2ds&=\displaystyle\int_{\partial \Omega_{+,\epsilon}}
e^{-x\cdot 2s\eta}|\partial_{\nu}(e^{-\frac{\phi_{1s}}{2}}u_0)|^2ds\\
&\lesssim \displaystyle\int_{\partial \Omega_{+,\epsilon}}
e^{-x\cdot 2s\eta}|\partial_{\nu}u_0|^2ds\\
&\lesssim \displaystyle\int_{\partial \Omega_{+,\epsilon}}
e^{-x\cdot 2s\eta}|\partial_{\nu}u|^2ds+\displaystyle\int_{\partial \Omega_{+,\epsilon}}
e^{-x\cdot 2s\eta}|\partial_{\nu}\delta u|^2ds,
\end{array}
\eqno{(5.14)}
\]
where we use the relation $u_0$\textbar$_{\partial \Omega}=0$ in the first inequality.

By Lemma $5.3$ and Lemma $5.4$ the inequality $(5.14)$ implies
\[
 \begin{array}{lll}
\lim\limits_{s\rightarrow\infty}\displaystyle\int_{\partial \Omega}\gamma_1 \partial_{\nu}(\tilde{u}_1-u_2)u_1ds=0.
\end{array}
\]

Then Lemma $5.5$ follows.

\bigbreak
\begin{center}
\textbf{\Large { The  Proof of Theorem $1.1$ }}
\end{center}

\textbf{Proof}\ \ Recalling $u_i^{(n)}=e^{-\frac{\phi_{is_n}}{2}}e^{x\cdot \zeta_i^{(n)}}(1+w_i^{(n)}), i=1,2.$ and using Proposition $5.1$
we know
 \[
 \begin{array}{lll}
\displaystyle\int_{\Omega}(\gamma_1^{\frac{1}{2}}\nabla \gamma_2^{\frac{1}{2}}-\gamma_2^{\frac{1}{2}}\nabla
\gamma_1^{\frac{1}{2}})\cdot \nabla (e^{\frac{\phi_{1s_n}}{2}}e^{x\cdot \zeta_1^{(n)}}(1+w_1^{(n)})e^{\frac{\phi_{2s_n}}{2}}e^{x\cdot \zeta_2^{(n)}}(1+w_2^{(n)}))dx\\
\ \ \ \ \ \ \ \ \ \ \ \ \ \ \ \  \ \  \ \  \ \ \ \ \ \  \ =\displaystyle\int_{\partial \Omega}\gamma_1 \partial_{\nu}(\tilde{u}_1^{(n)}-u_2^{(n)})u_1^{(n)}ds.
\end{array}
\eqno{(5.15)}
\]
Observing $e^{x\cdot\zeta_1^{(n)}}e^{x\cdot\zeta_2^{(n)}}=e^{ix\cdot k}$ and $\lim\limits_{n\rightarrow\infty}e^{-\frac{\phi_{is_n}}{2}}=\gamma_i^{-\frac{1}{2}}$,
by Lemma $5.5$ we have
 \[
 \begin{array}{lll}
&-\displaystyle\int_{\Omega}(\gamma_1^{\frac{1}{2}}\nabla \gamma_2^{\frac{1}{2}}-\gamma_2^{\frac{1}{2}}\nabla \gamma_1^{\frac{1}{2}})\cdot \nabla (\gamma_1^{\frac{1}{2}}\gamma_2^{\frac{1}{2}}e^{ix\cdot k})dx\\
&=\lim\limits_{n\rightarrow \infty}\displaystyle\int_{\Omega}(\gamma_1^{\frac{1}{2}}\nabla \gamma_2^{\frac{1}{2}}-\gamma_2^{\frac{1}{2}}\nabla \gamma_1^{\frac{1}{2}} )\nabla (e^{-\frac{\phi_{1s_n}}{2}}e^{-\frac{\phi_{2s_n}}{2}}e^{ix\cdot k})(w_1^{(n)}+w_2^{(n)}+w_1^{(n)}w_2^{(n)})dx\\
&+\lim\limits_{n\rightarrow \infty}\displaystyle\int_{\Omega}(\gamma_1^{\frac{1}{2}}\nabla \gamma_2^{\frac{1}{2}}-\gamma_2^{\frac{1}{2}}\nabla \gamma_1^{\frac{1}{2}}) (e^{-\frac{\phi_{1s_n}}{2}}e^{-\frac{\phi_{2s_n}}{2}}e^{ix\cdot k})(\nabla w_1^{(n)}+\nabla w_2^{(n)}+\nabla(w_1^{(n)}w_2^{(n)}))dx\\
&=:S_1+S_2.
\end{array}
\]

For $S_1$, recalling $\|\phi_{is_n}\|_{L^\infty}\leq \|\phi_{i}\|_{L^\infty}$ and the definition of $\Phi_{B}$ we have
 \[
 \begin{array}{lll}
|S_1|&\lesssim \limsup\limits_{n\rightarrow \infty}\displaystyle\int_{\mathbb{R}^n}(|\Phi_B w_1^{(n)}|+|\Phi_Bw_2^{(n)}|+|\Phi_Bw_1^{(n)}||\Phi_B w_2^{(n)}|)dx\\
&\lesssim \limsup\limits_{n\rightarrow \infty}(\|\Phi_B w_1^{(n)}\|_{L^2}+\|\Phi_B w_2^{(n)}\|_{L^2}+\|\Phi_B w_1^{(n)}\|_{L^2}\|\Phi_B w_2^{(n)}\|_{L^2}).
\end{array}
\eqno{(5.16)}
\]

On the other hand, by Lemma $3.2$ and Lemma $3.4$ we have
 \[
 \begin{array}{lll}
\|\Phi_B w_i^{(n)}\|_{L^2}\lesssim \frac{1}{s_n^{\frac{1}{2}}} \|\Phi_B w_i^{(n)}\|_{\dot{X}_{\zeta_i^{(n)}}^{\frac{1}{2}}}.
\end{array}
\]

 It follows that $S_1=0$.

For $S_2$, recalling $\gamma_1, \gamma_2\in H^{\frac{3}{2}, 2}(\Omega)$, $\|\phi_{is_n}\|_{L^\infty}\leq \|\phi_{i}\|_{L^\infty}$ and
$\|\nabla\phi_{is_n}\|_{L^\infty}\leq \|\nabla\phi_{i}\|_{L^\infty}$ we have
 \[
 \begin{array}{lll}
|S_2|&\lesssim \limsup\limits_{n\rightarrow \infty}(\|\Phi_B w_1^{(n)}\|_{H^{\frac{1}{2}, 2}}+\|\Phi_Bw_2^{(n)}\|_{H^{\frac{1}{2}, 2}}\\
&+\|\Phi_Bw_1^{(n)}\|_{L^2}\|\nabla(\Phi_B w_2^{(n)})\|_{L^2}+\|\nabla(\Phi_Bw_1^{(n)})\|_{L^2}\|\Phi_B w_2^{(n)}\|_{L^2}).
\end{array}
\eqno{(5.17)}
\]

Using Proposition $3.2$ and Lemma $3.3$ we have
 \[
 \begin{array}{lll}
&\|\Phi_B w_i^{(n)}\|_{H^{\frac{1}{2}, 2}}\lesssim \|w_i^{(n)}\|_{\dot{X}_{\zeta_i^{(n)}}}^{\frac{1}{2}},\\
&\|\Phi_B w_i^{(n)}\|_{L^2}\lesssim \frac{1}{s_n^{\frac{1}{2}}}\|w_i^{(n)}\|_{\dot{X}_{\zeta_i^{(n)}}}^{\frac{1}{2}},\\
&\|\nabla(\Phi_B w_i^{(n)})\|_{L^2}\lesssim s_n^{\frac{1}{2}} \|w_i^{(n)}\|_{\dot{X}_{\zeta_i^{(n)}}}^{\frac{1}{2}},
\end{array}
\]

by Lemma $3.4$, which implies  $S_2=0$.

Therefore, we have
\[
 \begin{array}{lll}
\displaystyle\int_{\Omega}(\gamma_1^{\frac{1}{2}}\nabla \gamma_2^{\frac{1}{2}}-\gamma_2^{\frac{1}{2}}\nabla \gamma_1^{\frac{1}{2}})\cdot \nabla (\gamma_1^{\frac{1}{2}}\gamma_2^{\frac{1}{2}}e^{ix\cdot k})dx=0.
\end{array}
\]

It follows that
\[
 \begin{array}{lll}
\displaystyle\int_{\Omega}e^{ix\cdot k}(-\frac{ik}{2}\cdot \nabla(\text{log}\gamma_1-\text{log}\gamma_2)+\frac{1}{4}((\nabla \text{log}\gamma_1)^2-(\nabla \text{log}\gamma_2)^2))dx=0
\end{array}
\eqno{(5.18)}
\]
for $k\bot \eta$. However, since the DN maps agree on ${\partial\Omega_{-,\epsilon}}(\eta)$ for a fixed constant $\epsilon>0$, they also 
agree on ${\partial\Omega_{-,\epsilon'}}(\eta')$ for $\eta'$ sufficiently close to $\eta$ on the unit sphere and for some smaller constant
 $\epsilon'$. Thus, in particular, $(5.18)$  holds for $k$ in an open cone in $\mathbb{R}^n$. Let the distribution $q$
be equal to $\frac{1}{2}\triangle(\text{log}\gamma_1-\text{log}\gamma_2)+\frac{1}{4}((\nabla\text{log}\gamma_1)^2-(\nabla\text{log}\gamma_2)^2)$ in $\Omega$
 and  zero outside of $\Omega$. Hence, $(5.18)$ implies that the 
Fourier transform of $q$ vanishes in an open set. Since $q$ is compact supported,
 the Fourier transform is analytic by the Paley-Wiener theorem, and this implies $q\equiv 0$. i.e.,
 \[
 \begin{array}{lll}
\frac{1}{2}\triangle(\text{log}\gamma_1-\text{log}\gamma_2)+\frac{1}{4}((\nabla\text{log}\gamma_1)^2-(\nabla\text{log}\gamma_2)^2)=0, \ \ \text{in}\ \Omega
\end{array}
\]
and since $\text{log}\gamma_1$\textbar$_{\partial \Omega}=\text{log}\gamma_2$\textbar$_{\partial \Omega}$, uniqueness of boundary value problem of
uniform elliptic equations implies $\gamma_1=\gamma_2$. Then the Theorem $1.1$ holds.

\textbf{Acknowledgment}\ \ The author wants to express his gratitude to Professor Mikko Salo who pointed out a mistake in the original version of this paper. 
The author also wants to thank the referees for many nice suggestions for the presentation of the paper.  This author was  supported by the Academy of Finland.

\bibliographystyle{alpha}
\bibliography{INVERSEPROBLEM}

\end{document}